\documentclass[12pt]{amsart}
\usepackage{amsmath,amsfonts,amssymb,amsthm,amscd,latexsym,euscript,xypic,
verbatim} \hoffset=-1cm \textwidth=15.5cm \theoremstyle{plain}

\newtheorem*{T1}{Theorem 1}
\newtheorem*{T2}{Theorem 2}

\newtheorem*{thm1.2}{(1.2) Theorem}
\newtheorem*{thm1.3}{(1.3) Theorem}
\newtheorem*{thm1.4}{(1.4) Theorem}
\newtheorem*{propA*}{Proposition A}
\newtheorem*{propB*}{Proposition B}
\newtheorem*{thmC*}{Theorem C}
\newtheorem*{propD*}{Proposition D}
\newtheorem{prop}{Proposition}[section]
\newtheorem{thm}[prop]{Theorem}
\newtheorem{cor}[prop]{Corollary}
\newtheorem{lemma}[prop]{Lemma}

\theoremstyle{definition}

\newtheorem{Def}[prop]{Definition}
\newtheorem*{Def*}{Definition}
\newtheorem{Defs}[prop]{Definitions}

\newtheorem{notation}[prop]{Notation}
\newtheorem*{notation*}{Notation}
\newtheorem*{question*}{Question}
\newtheorem{remark}[prop]{Remark}
\newtheorem{remarks}[prop]{Remarks}
\newtheorem*{rem}{Remark}

\newcommand{\calb}{\mathcal{B}}

\newcommand{\qq}{\mathbb{Q}}

\newcommand{\zz}{\mathbb{Z}}

\newcommand{\ga}{\alpha}
\newcommand{\gb}{\beta}

\newcommand{\gl}{\lambda}

\newcommand{\gvp}{\varphi}

\newcommand{\sminus}{\smallsetminus}
\newcommand{\lan}{\langle}
\newcommand{\ran}{\rangle}

\newcommand{\End}{{\rm End}}

\newcommand{\Hom}{{\rm Hom}}

\newcommand{\ch}{\check}

\newcommand{\onto}{\twoheadrightarrow}

\newcommand{\longright}[1]{\;{\count255=0 \loop \relbar\mathrel{\mkern-6mu}%
   \advance\count255 by1\ifnum\count255<#1\repeat\rightarrow}\;}
\numberwithin{equation}{section}
\begin{document}
\title[On Kernels of cellular covers]{On kernels of cellular covers}
\author[Farjoun, G\"{o}bel, Segev, Shelah]
{E.~D.~Farjoun$^1$\quad R.~G\"{o}bel$^2$\quad Y.~Segev$^3$\quad
S.~Shelah$^4$}

\address{Emmanuel D.~Farjoun\\
        Institute of Mathematics\\
        The Hebrew University of Jerusalem, Givat Ram\\
        Jerusalem 91904\\
        Israel}
\email{farjoun@math.huji.ac.il}

\address{R\"{u}diger G\"{o}bel\\
        Fachbereich Mathematik\\
        Universit\"{a}t-Duisburg-Essen\\
        D-45117 Essen\\
        Germany}
\email{R.Goebel@uni-essen.de}

\address{Yoav Segev\\
        Department of Mathematics\\
        Ben-Gurion University\\
        Beer-Sheva 84105\\
        Israel}
\email{yoavs@math.bgu.ac.il}

\address{Saharon Shelah\\
Institute of Mathematics\\
        The Hebrew University of Jerusalem, Givat Ram\\
        Jerusalem 91904\\
        Israel 
        and Department of Mathematics \\
        Rutgers University\\ 
        New Brunswick, NJ 08854-8019, USA} 
\email{shelah@math.huji.ac.il}

%
%
%
%
\begin{abstract}
In the present paper we continue to examine cellular covers
of groups, focusing on the cardinality and the structure
of the kernel $K$ of the cellular map $G\to M$.  We 
show that in general a  torsion free reduced abelian group $M$ may have
a proper class of non-isomorphic cellular covers. In other words,
the cardinality of the kernels is unbounded. In the opposite direction
we  show
that if the kernel of a cellular cover of any group $M$
has certain ``freeness'' properties, then its cardinality
must be bounded by $|M|$.
\end{abstract}
\date{November 23, 2006}
\dedicatory{Dedicated to Avinoam Mann on the occasion of his retirement, 2006.}
\maketitle

\footnotetext[1]{Partially supported by a grant from
the Israel Science Foundation.} 
\footnotetext[2]{Partially
supported by project no.~I-706-54.6/2001 of the German-Israeli
Foundation for Scientific Research \& Development.}
\footnotetext[3]{Partially supported by BSF grant no.~2004-083.}
\footnotetext[4]{Partially
supported by project no.~I-706-54.6/2001 of the German-Israeli
Foundation for Scientific Research \& Development; this is the author's
paper \# 892.}
\renewcommand{\thefootnote}{}

\renewcommand{\thefootnote}{\arabic{footnote}}
\vspace{1.3cm}
%
%
%
%
%

\section*{Introduction and main results}\label{S:ii}

In this paper we continue the discussion of \text{cellular covers}
in the category of groups begun in \cite{FGS1,FGS2}, where this
notion is also motivated. Given a map of groups $c\colon  G\to M,$
we say that $(G,c)$ is {\em a cellular cover of $M$} or that
$c\colon G\to M$ is a cellular cover, if every  group map
$\gvp\colon G\to M$ factors uniquely through $c$, or,
equivalently, the natural map $\Hom(G,G)\to \Hom(G,M),$ induced by
$c,$  is an isomorphism of sets. Explicitly this means that there
exists a unique \textit{lift} $\widetilde\gvp\in\End(G)$ such that
$\widetilde\gvp\circ c=\gvp$ (maps are composed from left to
right).

It has been shown before \cite{FGS1,FlR} that cellular covers are
values of general augmented ($FM\to M$) and idempotent
($F\circ F=F$) functors on the category of groups. More
concretely, such functors are of the form $cell_A(-),$ namely
$A$-cellular approximation with respect to some group $A$.

The functors $cell_A(-)$ had been used fruitfully in the category
of groups, topological spaces and chain complexes over rings or
DGAs ($=$ differential graded algebras); compare, for example,
Dwyer {\it et al} \cite{DGrI},  \cite{RSc,FlR}, Shoham
(see \cite{Sho}). The present results shade some light on the possible
values of the functor $cell_A(-)$ when $A$ is abelian. (We note
that very different groups $A$ can give rise to the same
functor.)  
 It is possible that the values of all such functors 
 (i.e.~$\{cell_AM\mid A\text{ a group}\}$)
 on a fixed group $M$ yields only a set of results, up to isomorphism. In
some topological analogous situations it has been shown that
indeed only a set of values occurs (cf.~\cite{DP}). We have seen in
\cite{FGS1,FGS2} that this is the case when $M$ 
is a finite group, a finitely generated nilpotent group or
a divisible abelian group. One aim of the present paper is to show
that {\em there are abelian groups $M$ for which 
$\{cell_AM\mid A\text{ an abelian group}\}$
is a proper
class of
isomorphism types}.  This is a consequence of the following.  
\begin{T1}\label{T1}
For any infinite cardinality $\gl$,  there exists an abelian  
group $M$ of cardinality $\gl$ with $\End\, M\cong \zz,$ such
that for any infinite cardinality $\kappa$ there exists an abelian
group $K$ of cardinality $\kappa$ and with $\Hom(K,M)=0$ such
that $K$ is the kernel of some cellular cover $G\onto M.$
\end{T1}
\noindent
Theorem 1 is Theorem \ref{cc.11} of \S2;  its proof
relies on Theorem \ref{cc.5} which may be of independent
interest.

Let $c\colon G\to M$ be a cellular cover.
In previous papers we have noticed that $G$ inherits several
important properties from $M$: First the kernel $K=\ker c$ is
central in $G$, that is, $G$ is a central extension of $M,$ and
further, if $M$ is nilpotent, then $G$ is nilpotent of the same
class; if $M$ is finite then so is $G$. In addition, we have
classified all possible covers of divisible abelian groups
(\cite[\S 4]{FGS2}) and showed that when $M$ is abelian the
kernel $K$ is reduced and torsion-free (\cite[Thm. 4.7]{FGS1}.
The case when $M$ is abelian was independently investigated in
\cite{BD} and \cite{D}.  Amongst other results it was shown there
that when $M$ is (abelian and) reduced, $K$ is cotorsion free.

In \cite{FGS1} we have already observed that if $M$ is perfect,
and $G$ is the so-called universal central extension of $M$ (so
that $K$ is the Schur-multiplier), then $G\twoheadrightarrow M$ is
a cellular cover, and, since {\it any} abelian group is a Schur-multiplier, 
in general, there is no restriction on the {\it
structure} of $K$ (other than being in the center of $G$ and hence
$K$ is abelian).

Note that the covers in Theorem 1
are very special covers in which the  only map $K\to M$ from the
kernel to $M$ is the zero map. This class of maps are both
cellular cover and localization maps. Namely $c\colon G\to M$ is
both a cellular cover and a localization.  Recall that ``$c$ is a
localization'' means that for any $\gvp\in\Hom(G,M)$ there is a
unique {\it corresponding} $\widetilde\gvp\in\End(M)$ such that
$c\circ\widetilde\gvp=\gvp$. Therefore, this class of
localization-cellular maps $V\to W$ have the property that they
induce isomorphisms on endomorphism sets: $\End\,V\cong \Hom
(V,W)\cong \End\,W$.

 The kernel $K$ in Theorem 1 cannot be of an arbitrary nature:
\begin{T2}\label{T2}
For any cellular cover $c:G\to M$ (where $M$ is an arbitrary,
not necessarily abelian,
group), if the kernel  $K$ of $c$ is a free abelian group then
$|K|\le |M|.$
\end{T2}
In fact, the results in \S 1 (see Proposition \ref{bb.4}) are somewhat more general than Theorem 2.  We note that \cite{FuG} continues the investigation
of cellular covers of abelian groups begun in \cite{FGS2}
and in Theorem 1 of this paper, and in particular, further results
on ``large'' cellular covers
of ``small'' abelian groups are obtained there.  
%
%
%
%
%
\section{Free kernels  are small}\label{S:bb}
In this section we consider the  kernel $K$ of a cellular cover
$c\colon G\to M$. We impose some additional
``freeness'' assumptions on $K$. We show that under these
restrictions the cardinality of $G$ is bounded in terms of the
cardinality of $M$.
\begin{Def}[Compare with \cite{EMe}, p.~90, \cite{Fu}, p.~184]\label{bb.1}
Let $K$ be an abelian group and $\ga, \gb$ be cardinal numbers
such that $\ga\le\gb$. We say that $K$ is \textit{weakly-$(\ga,
\gb)$-separable} iff any subgroup $K_1\le K$ of size $\le \ga$ is
contained in a direct summand $K_2\le K$ of size $\le \gb$.
Notice that when $\ga=\gb$, then our notion coincides with the
notion of (weakly) $\ga^+$-separable group as in \cite{EMe}, p.~90.
In this case we will say that $K$ is weakly-$\ga$-separable (and
not weakly $\ga^+$-separable as in \cite{EMe}).
\end{Def}
We recall the following well-known fact.
\begin{lemma}\label{bb.2}
Let $K$ be a free abelian  group.  Then $K$ is
weakly-$\ga$-separable, for every infinite cardinal number $\ga$.
\end{lemma}
\begin{proof}
Let $K_1$ be a subgroup of $K$.  Of course we may assume that
$K_1\ne 0$. Let $\calb$ be a basis of $K$ and for each $x\in K_1$
let $\calb_x\subseteq\calb$ be a finite subset such that
$x\in\lan\calb_x\ran$.  Let $K_2:=\lan \calb_x\mid x\in K_1\ran$.
Then $K_1\le K_2$, $\mid K_1\mid =\mid K_2\mid$, and $K=K_2\oplus
F$, where $F=\lan \calb\sminus\bigcup_{x\in K_1}\calb_x\ran$.
\end{proof}
\begin{lemma}\label{bb.3}
If $G, M$ are groups and $c\in\Hom(G,M)$ is surjective, then there
exists $G_1\le G$ such that $\mid G_1\mid \le \mid M\mid+\aleph_0$
with $c(G_1)=M$.
\end{lemma}
\begin{proof}
For each $m\in M$ choose a preimage $g_m\in G$ (i.e.~$c(g_m)=m$)
and let $G_1=\lan g_m\mid m\in M\ran$.
\end{proof}
\begin{prop}\label{bb.4}
Let $c\colon G\to M$ be a cellular cover of the infinite group $M$
and set $K:=\ker c$. Let $\gb$ be a cardinal number such that
$\gb\ge\mid M\mid$.  Then
\begin{itemize}
  \item[(1)]  if $K$ is weakly $(\mid M\mid, \gb)$-separable, then
  $\mid G\mid\le \gb$; in particular,
 \item[(2)]  if $K$ is a free abelian group, then $\mid G\mid\le\mid
  M\mid$.
\end{itemize}
\end{prop}
\begin{proof}
Notice that (2) is an immediate consequence of (1) and Lemma
\ref{bb.2}. It remains to prove (1). 
Notice that if we restrict the image and consider the
map $c\colon G\to c(G)$ we still get a cellular cover.  It follows that
if $c(G)$ is finite, then $G$ is finite (see \cite[Theorem 5.4]{FGS1}).
We may thus assume without loss
that $c$ is surjective. Let $G_1\le G$ be a subgroup such that
$c(G_1)=M$ and such that
\[
\mid G_1\mid =\mid M\mid,
\]
whose existence is guaranteed by Lemma \ref{bb.3} (note that since
$M$ is infinite, $|M|+\aleph_0=|M|$). Since $c(G_1)=M$, we have
that
\[
G=KG_1.
\]
Let $K_1:=G_1\cap K;$ then $\mid K_1\mid\le\mid M\mid$, so by
hypothesis there exists a subgroup $K_2\le K$ such that $K_1\le
K_2$, $\mid K_2\mid\le\gb$ and such that $K=K_2\times F$, for some
$F\le K$. It is easy to check that it follows that
\[
G=(G_1K_2)\times F.
\]
In particular, if $F\ne 1$, then, since $F\le K$, $\Hom(G,K)\ne 0,$
a contradiction. Thus $F=1$, so $G=G_1K_2$ and hence $\mid
G\mid\le\gb$.
\end{proof}
%
%
%
%
\section{Cellular covers  with large kernels}\label{S:cc}
\medskip

\subsection*{A. Preliminaries}\hfill
\medskip

Before  describing the main construction  we  introduce some
definitions, prove a few lemmas about them  and   recall an
existence result about ``large'' rigid abelian groups to be used below.
\begin{Defs}\label{cc.1}
Let $A$ be an abelian group, $q$ a prime and $\pi$ a set of primes. Then
\begin{enumerate}
    \item $A$ is {\it $q$-reduced} if
$\bigcap_{i=1}^{\infty}q^iA=0$.
   \item $A$ is $\pi$-reduced if $A$ is $p$-reduced, for all $p\in \pi$.
    \item An element $a\in A$ is $q$-pure (in $A$) if $a$ is
not divisible by $q$ in $A$.
    \item $A$ is $q$-divisible if each element $a\in A$ is
divisible by $q$ in $L$.

   \item An integer $n$ is a $\pi$-number, if $n$ is divisible only by
   primes from $\pi$ ($1$ and $-1$ are always $\pi$-numbers).
   
   \item A torsion element $a\in A$ is a $\pi$-element if the order of $a$ is a $\pi$-number
   (or $a=0$).
   \item $A$ is a $\pi$-group, if each element of $A$ is a $\pi$-element. 
  \item $\zz[1/\pi]:=\zz[1/p\mid p\in \pi]$ (and if $\pi=\emptyset$, then $\zz[1/\pi]=\zz$). 
\end{enumerate}
\end{Defs}
\begin{remarks}[Tensor products, see \cite{Fu}]\label{cc.2}
\begin{itemize}
\item[(1)]
Let $A$ be a torsion free abelian group.  Then $V:=\qq\otimes A$ is a
vector space over $\qq$ which contains a copy of $A$.  Thus we always
think of $A$ as being contained in a vector space $V$ over $\qq$
such that $V/A$ is a torsion abelian group.
Hence it makes sense to talk about
the group $\lan A\cup \{\frac{a_i}{m_i}\mid i\in I\}\ran$
where $I$ is an index set, $\{a_i\mid i\in I\}\subseteq A$ and
$\{m_i\mid i\in I\}\subseteq \zz\sminus\{0\}$.  This is the subgroup
of $V$ generated by $A\cup \{\frac{a_i}{m_i}\mid i\in I\}$.
\item[(2)]
Note that if $S\subseteq V$
and $\pi$ is a set of primes  such that for each $s\in S$ there exists
a $\pi$-number $n$ with $ns\in A$, then $\lan A\cup S\ran/A$ is
a $\pi$-group.  In particular, for a subring $R\subseteq Q$ we view
$R\otimes A$ as a subgroup of $V$ and if $R=\zz[1/\pi]$, then $(R\otimes A)/A$
is a $\pi$-group.
\item[(3)]
Note further that if $\pi_1$ and $\pi_2$ are disjoint sets
of primes  and $B\subseteq V$ is a subgroup containing $A$  such that
$A$ is $\pi_1$-reduced and $B/A$ is a $\pi_2$-group, then $B$
is $\pi_1$ reduced.   
\end{itemize}
\end{remarks}
\begin{notation}\label{cc.3}
Let $L$ be a torsion free abelian group and let $q$ be a prime.
Let $0\ne x\in L$ we denote, using Remark \ref{cc.2}(1),
\[
\textstyle{L\oplus_x\zz[1/q]:=\lan L\cup \{\frac{x}{q^i}\mid 1\le i\in\zz\}\ran}.
\]
We write $H=x\zz[1/q]$ for the subgroup of  
$L\oplus_x\zz[1/q]$ consisting of the elements
\[
\textstyle{H:=\{\frac{m}{q^i}x\mid m\in\zz\text{ and }1\le i\in\zz\}.}
\]
\end{notation}
\begin{remark}\label{cc.4}
Assume $L$ is a torsion free abelian group, $q$ is a prime
and $0\ne x\in L$ is a $q$-pure element.  Then
\[
L\oplus_x\zz[1/q]\cong (L\oplus\zz[1/q])/\lan (-x, 1)\ran.
\]
Furthermore,
let $\widehat M$ be a group such that $\widehat M=L\oplus H$ where
$L, H$ are subgroups of $\widehat M$, $L$ is torsion free and $H$
is isomorphic to $\zz[1/q]$ under an isomorphism taking some $0\ne h\in H$ to
$1$.  Let $0\ne y\in L$ be a $q$-pure element and let $M:=\widehat M/\lan y-h\ran$. Then
$M$ is isomorphic to the group $L\oplus_y\zz[1/q]$ constructed in
Notation \ref{cc.3}.
\end{remark}
\medskip

\subsection*{B. Existence of large rigid groups}\hfill
\medskip

The following is our main stepping stone for proving the existence
of covers with arbitrarily large kernels.

\begin{thm}\label{cc.5}
Let $P$ be a set of at least  four primes, $Q$ its complementary
set of primes and $\gl$ any infinite cardinal. Then there is a
torsion-free abelian group $H$ of cardinality $\gl$ with the
following three properties.
\begin{enumerate}
 \item $H$ is $Q$-reduced;
 \item  if $Q_0\subseteq Q$ is a set of primes and 
 $A$ is a torsion free  abelian group containing $H$
 such that $A/H$ is a $Q_0$-group, then 
  $\End(A)\subseteq \zz[1/Q_0]$;
  \item  $H$ contains a free abelian group $F$ of cardinality
  $\gl$ such that $H/F$ is a $P$-group.
\end{enumerate}
\end{thm}
\begin{proof}
Let $R:=\zz[1/Q]$.
By \cite[Thm.~2.1]{Sh} (see also \cite[Corollary 14.5.3(b), p.~577]{GT}),
there exists an $R$-module $M$ of cardinality $\gl$ such that
$\End(M)=R$.  Let $\calb$ be a maximal $(\zz)$-independent
subset of $M$.  We let
\[
F:=\lan \calb\ran\text{ and } H:=\{x\in M\mid \text{there exists a $P$-number $n\in\zz$
with $nx\in F$}\}.
\]
We claim that $H$ satisfies all the required properties.
By construction (3) holds.  Also, since $F$ is a free abelian group
and since $H/F$ is a $P$-group, $H$ is $Q$-reduced (see Remark \ref{cc.2}(3)),
so (1) holds.  

We now show (2).  By construction, $M/H$ is a $Q$-group,
so $R\otimes H=M$.  Thus for any group $H\subseteq A\subseteq M$, $R\otimes A=M$.
Let $A$ be as in (2).  Then $H\subseteq A\subseteq R\otimes A=M$, and
since $R\otimes A=M$, it follows that any endomorphism of $A$ extends to
an endomorphism of $M$, thus $\End(A)\subseteq R$.  Let
now $Q_0\subseteq Q$ and suppose that $A/H$ is a $Q_0$-group.
Let $f\in\End(A)$ so that $f$ is multiplication by $\frac{m}{n}$,
where ${\rm gcd}(m,n)=1$ and $n$ is a $Q$-number.  Assume 
there exists a prime $q\in Q\sminus Q_0$ such that
$q\mid n$.  Then, after multiplying by an appropriate integer,
we may assume that $n=q$.  Writing $1=\ga q+\gb m$, with
$\ga, \gb\in\zz$, we see that $\frac{1}{q}=\ga+\gb\frac{m}{q}$,
so multiplication by $\frac{1}{q}$ is an endomorphism of $A$.
However, $q\notin Q_0$, $H$ is $Q$-reduced and $A/H$ is a $Q_0$-group,
so Remark \ref{cc.2}(3) implies that $A$ is $q$-reduced.  This is a contradiction.
Thus $n$ is a $Q_0$-number, so $\End(A)\subseteq \zz[1/Q_0]$ and (2) holds.

 \end{proof}
\begin{rem}
The set of primes $P$ in Theorem \ref{cc.5} is
the set of primes that are used to construct the
$\zz[1/Q]$-module $M$ as in the begining of the proof
of the theorem.  Thus we only work with the complimentary set
of primes $Q$ when using the theorem to construct groups
$L$ that have some desirable properties.  Below we
will fix the set $Q$ of primes which will be used 
for our constructions (in fact we only need $3$ primes in $Q$, see Corollary \ref{cc.6} below).  The set $P$ will be the complimentary set of primes.
\end{rem}
The variant of Theorem \ref{cc.5} which we actually use in
subsection C below is the following Corollary.
\begin{cor}\label{cc.6}
Let $\gl$ be any infinite cardinal and let $Q:=\{q_L, q_K, q\}$ be
a set consisting of three primes.
Then there exists an abelian group $L$ whose
cardinality is $\gl$ such that
\begin{itemize}
 \item[(1)]  $L$ is  torsion free and $q_L$-divisible; 
 \item[(2)]  $L$ is $Q\sminus\{q_L\}$-reduced;
 \item[(3)] if $M\supseteq L$ is a torsion free abelian group
 such that $M/L$ is a $q$-group, then $\End(M)\subseteq \zz[1/\{q_L, q\}]$.
 \item[(4)] there exists a $q$-pure element $x_L\in L$ such
 that for $M:=L\oplus_{x_L} \zz[1/q]$ we have $\bigcap_{i=1}^{\infty}q^iM=x\zz[1/q]$.
\end{itemize}
 \end{cor}
\begin{proof}
We use Theorem \ref{cc.5} with $Q$  playing the
role of $Q$ in that theorem.  
Let $H$ be as in Theorem \ref{cc.5},
let $R=\zz[1/q_L]$ and let $L:=R\otimes H$.
Notice that by Remark \ref{cc.2}(3), $L$ is $Q\sminus\{q_L\}$ reduced.
Of course $L$ is $q_L$-divisible.

Next if $M\supseteq L$ is a torsion free abelian group
such that $M/L$ is a $q$-group, then, by construction,
$M/H$ is a $\{q_L, q\}$-group, so (3) follows from
Theorem \ref{cc.5}(2).

To prove (4)  let $F$ be as in part (3) of Theorem \ref{cc.5}. Let 
$\calb\subseteq F$ be a  free generating set of $F$, pick $x_L\in\calb$
and set $x:=X_L$.  Clearly $x$ is $q$-pure.  Assume (4) is false and write
$U:=\bigcap_{i=1}^{\infty} q^i(L\oplus_x\zz[1/q])$.
Since $x\zz[1/q]\subseteq U$, there exists 
$\ell\in L\sminus\lan x\ran$ such that  $\ell\in U$.  But then
writing $\ell=\sum_{i=1}^t\ga_ix_i$, with $\ga_i\in\zz[1/(P\cup \{q_L\})]$, 
$x_i\in \calb$ and
$x_1\ne x$, we see that there exists $0 < j\in \zz$ such that $q^j$ does
not divide $\ga_1x_1$, and hence $q^j$ does not divide $\ell+sx$, for any
$s\in\zz$, and this contradicts the fact that $\ell\in U$.
\end{proof}
\medskip

\subsection*{C. Constructing  covers with arbitrarily large kernels}\hfill
\medskip

In this section we use Corollary \ref{cc.6} above  to construct an abelian
group $M$ and, for arbitrarily large cardinal $\kappa$, a cellular cover $G\to
M $ whose kernel $K$ has  cardinality $\kappa$. 
The group $M$ will be as in Corollary \ref{cc.6}(4).
Lemma \ref{cc.8} below describes the nice properties
of such a group $M$.

We start with a very simple lemma that allows us to conclude that
the canonical homomorphism $G\to G/K$ from the abelian group $G$ to
the factor group $G/K$
is a cellular cover. The rest of the
section is devoted to building arbitrarily large groups $K$ satisfying
the conditions of the lemma (while $G/K$ remains fixed).
\begin{lemma}\label{cc.7}
Let $G$ be an abelian group and $K\le G$ be a subgroup.  Set
$M:=G/K$ and let $c\colon G\to G/K$ be the canonical homomorphism.
Assume that
\begin{itemize}
   \item[(i)]  $\End(M)\cong\zz$;
   \item[(ii)]  $K$ is a fully invariant subgroup of $G$;
  \item[(iii)]  $\Hom(K,M)=0=\Hom(G,K)$.
\end{itemize}
Then $\End(G)=\zz$ and $c$ is a cellular cover.
\end{lemma}
\begin{proof}
Let $\mu\in\End(G)$.  By (ii), $\mu(K)\le K$ so $\mu$ induces
$\hat\mu\in\End(M)$ defined by $\hat\mu(g+K)=\mu(g)+K$.  By (i),
there exists $n\in\zz$ such that $\hat\mu$ is multiplication by
$n$.  Thus the map $g\to (\mu(g)-ng)$ is in $\Hom(G,K)$, so by
(iii) it is the zero map and it follows that $\mu$ is
multiplication by $n$.  This shows that $\End(G)\cong\zz$.

Let now $\gvp\in\Hom(G,M)$.  Then by (iii), $\gvp(K)=0$, so $\gvp$
induces $\hat\gvp\in\End(M)$ defined by $\hat\gvp(g+K)=\gvp(g)$.
Thus by (i) there is $n\in\zz$ such that $\gvp(g)=ng+K$, for all
$g\in G$. Consequently, the map $\ch\gvp\in\End(G)$ defined by
$\ch\gvp(g)=ng$ lifts $\gvp$, so any $\gvp\in\Hom(G,M)$ lifts.
Since $\Hom(G,K)=0$,  \cite[Lemma 3.6]{FGS1} shows that $c$ is a
cellular cover.
\end{proof}
%
%
%
%
%
\begin{lemma}\label{cc.8}
Let $Q:=\{q_L, q_K, q\}$ be
a set consisting of three primes, and let $L$ be an abelian group
satisfying (1)--(4) of Corollary \ref{cc.6}. 
Let $x_L\in L$ be a $q$-pure element as in (4) of Corollary
\ref{cc.6}, and set $M=L\oplus_{x_L}\zz[1/q]$.  Then $M$ is torsion
free, it is $q_K$-reduced and $\End(M)\cong\zz$.
\end{lemma}
\begin{proof}
That $M$ is torsion free is by construction.
By  Remark \ref{cc.2}(3), $M$ is $q_K$-reduced.

Recall that by (4) of Corollary \ref{cc.6}, 
\begin{equation*}\tag{*}
H=\bigcap_{i=1}^{\infty}q^iM,
\end{equation*}
where $H=x\zz[1/q]$ is as in Notation \ref{cc.3}.

Let  $\gvp\in\End(M)$.
Since  $M/L$ is a $q$-group part (3) of Corollary \ref{cc.6}
implies that 
there exists $\frac{m}{n}\in\qq$, with ${\rm gcd}(m,n)=1$
such that $n\ge 1$ is a $\{q_L, q\}$-number and such that $\gvp(x)=\frac{m}{n}x$,
for all $x\in M$.  Suppose $n\ne 1$ and let $p\in\{q_L, q\}$ such that $p\mid n$.
Since $\frac{m}{n}x\in M$, for  all $x\in M$ also $\frac{m}{p}x\in M$,
for all $x\in M$ and then writing $1=\ga m+\gb p$, $\ga, \gb\in\zz$ we see
that $\frac{1}{p}x=\frac{\ga m }{p}x+\gb x\in M$.  Thus $M$ is $p$-divisible.
Now if $p=q$, then (*) implies that $L$ is not $q$-divisible,
a contradiction.
If $p=q_L$, then, since by (*) $H$ is a fully invariant subgroup of $M$,
it follows that $H$ is $q_L$-divisible (because multiplication
by $1/q_L$ is an endomorphism of $M$).  But of course $H$
is not $q_L$ divisible.  Thus $n=1$ and this completes the proof of the lemma.
\end{proof}
\begin{lemma}\label{cc.9}
Let $G$ be an abelian group containing  subgroups $K$ and
$\widehat M$  such that $G=K+\widehat M$.  Set $M:=G/K$ and let
$c\colon G\to M$ be the canonical homomorphism.  Assume that
\begin{itemize}
 \item[(i)]  $K$ is a torsion free fully invariant subgroup of $G$;
 \item[(ii)]  $K$ is an $R$-module for some subring $R\subset\qq$ and
$\End(K)=R$;
 \item[(iii)]  $M$ is torsion free and $\End(M)\cong\zz$;
 \item[(iv)]  $\Hom(\widehat M, K)=0=\Hom(K,M)$;
 \item[(v)]  $K\cap\widehat M\ne 0$.
\end{itemize}
Then $\End(G)\cong\zz$ and $c$ is a cellular cover.
\end{lemma}
\begin{proof}
We use Lemma \ref{cc.7}.  It only remains to show that
$\Hom(G,K)=0$. Let $\mu\in\Hom(G,K)$.  By hypothesis (iv),
$\mu(\widehat M)=0$. By hypothesis (i), $\mu(K)\le K$, so by
hypothesis (ii) there exists $r\in R$ such that $\mu(v)=rv$, for
all $v\in K$. Let $0\ne v\in\widehat M\cap K$.  Then
$rv=\mu(v)=0$, so since $K$ is torsion free, $r=0$, and it follows
that $\mu(K)=0$ and then $\mu=0$.
\end{proof}
\begin{prop}\label{cc.10}
Let $Q:=\{q_L, q_K, q\}$ be a set consisting of three primes.
Let $K$ and $L$ be abelian groups and assume that
\begin{itemize}
 \item[(i)] $K$ is torsion free, it is $q_K$-divisible
 and $Q\sminus \{q_K\}$-reduced.  
 \item[(ii)]  $L$ and the element $x_L\in L$ satisfy (1)--(4)
 of Corollary \ref{cc.6}   
\end{itemize}
Let  $0\ne x_K\in K$ be an
arbitrary element, and   let
\[
G=(K\oplus L)\oplus_{(x_K-x_L)}\zz[1/q]
\]
be the group constructed in Notation \ref{cc.3}, with $K\oplus L$
in place of $L$ and $x_K-x_L$ in place of $x$.
 Set
\[
H:=(x_K-x_L)\zz[1/q],\quad\text{ and }\quad\widehat M=L+H.
\]
Then $G$, $K$ and $\widehat M$ satisfy all the hypotheses of Lemma
\ref{cc.9}.  In particular, the canonical homomorphism $c\colon
G\to G/K$ is a cellular cover.
\end{prop}
\begin{proof}
Clearly $G=K +\widehat M$. 
Now since $(K+L)\cap H=\lan
x_K-x_L\ran$, it is easy to check that
\[\tag{I}
K\cap\widehat M=\lan x_K\ran.
\]
Note that $L\cap H=0$, because if $g:=n(x_K-x_L)/q^i\in L$, then
$n(x_K-x_L)\in L$, which implies that $nx_K\in L$.  But $K$ is
torsion free and $K\cap L=0$, so $n=0$ and then $g=0$. Thus
$\widehat M=L\oplus H$, also $x_K=x_L +(x_K-x_L)$ and
$H\cong\zz[1/q]$ by an isomorphism sending $(x_K-x_L) \to 1$, so
by (I) and Remark \ref{cc.4}, $M\cong \widehat M/\lan x_K\ran\cong
L\oplus_{x_L}\zz[1/q]$.
\relax From (ii)  and Lemma \ref{cc.8} it follows that
\[\tag{II}
M\text{ is torsion free, }M\text{ is $q_K$-reduced and
}\End(M)=\zz.
\]
Since $K$ is $q_K$-divisible, we conclude that
\[\tag{III}
\Hom(K,M)=0,
\]
and also, since $M$ is $q_K$-reduced, we have: $\bigcap_{i=0}^{\infty}q_K^iG=K$, so
\[\tag{IV}
\text{$K$ is a fully invariant subgroup of $G$}.
\]

Next, since $L$ is $q_L$-divisible and $K$ is $q_L$-reduced,
$\Hom(L,K)=0$.  Similarly, since $H$ is $q$-divisible,
$\Hom(H,K)=0$. Hence
\[\tag{V}
\Hom(\widehat M, K)=0.
\]
Thus all hypotheses of Lemma \ref{cc.9} have been verified.
\end{proof}
As a Corollary to Proposition \ref{cc.10} we get Theorem 1 of the
introduction.
\begin{thm}\label{cc.11}
Let $\gl$ be any infinite cardinal.  There exists an abelian group
$M$ of cardinality $\gl$ such that for any infinite cardinal
$\kappa\ge\gl$ there exists a cellular cover $c\colon G\to M$ with
$\mid \ker c\mid=\kappa$.
\end{thm}
\begin{proof}
Corollary \ref{cc.6} guarantees the existence of groups $L$ and
$K$ of cardinality $\gl$ and $\kappa$ respectively, and primes $q_K$,
$q_L$ and $q$ satisfying all hypotheses of Proposition
\ref{cc.10}. Let $K$ and $G$ be as in Proposition \ref{cc.10} and
set $M:=G/K$. By Proposition \ref{cc.10}, $c\colon G\to M$ is a
cellular cover and of course $\mid K\mid=\kappa$ and $\mid
M\mid=\gl$.  Notice that we saw in the proof of Proposition
\ref{cc.10} that $M\cong L\oplus_{x_L}\zz[1/q]$, so the structure
of $M$ is independent of the choice of $K$.
\end{proof}

\end{document}